\documentclass{amsart}
\usepackage{amssymb}
\usepackage[centertags]{amsmath}
\usepackage{amssymb}
\usepackage{xr}
\usepackage{amscd}
\usepackage{amsthm}
\usepackage{amsxtra}

\copyrightinfo{2004}{American Mathematical Society}

\newcommand{\re}{\ensuremath{\text{Re }}}
\newcommand{\im}{\ensuremath{\textrm{Im }}}

\newtheorem{theorem}{Theorem}[section]
\newtheorem{lemma}[theorem]{Lemma}

\theoremstyle{definition}
\newtheorem{definition}[theorem]{Definition}

\theoremstyle{remark}
\newtheorem{remark}[theorem]{Remark}

\numberwithin{equation}{section}

\begin{document}
\author{Albert Boggess}
\address{Department of Mathematics\\ Texas A \& M University\\
College Station Texas, 77843-3368} \email{boggess@math.tamu.edu}

\author{Daniel Jupiter}
\address{Department of Mathematics\\ Texas A \& M University\\
College Station Texas, 77843-3368} \email{jupiter@math.tamu.edu}

\subjclass[2000]{Primary 32V10, 32V99, 30E10}

\date{September 7, 2004}

\keywords{CR approximation, Bloom-Graham model graphs}

\commby{??}
\title[Global approximation of CR functions]{Global approximation of
CR functions on Bloom-Graham model graphs in $\mathbb{C}^n$}

\begin{abstract}
We define a class of generic CR submanifolds of $\mathbb{C}^n$ of
real codimension $d$, $1\leq d\leq n-1$, called the Bloom-Graham
model graphs, whose graphing functions are partially decoupled in
their dependence on the variables in the real directions. We prove
a global version of the Baouendi-Treves CR approximation theorem,
for Bloom-Graham model graphs with a polynomial growth assumption
on their graphing functions.
\end{abstract}

\maketitle

\section{Introduction}
Baouendi and Treves \cite{baouendi-treves:CR} proved that CR
functions on a generic CR submanifold of $\mathbb{C}^n$ can be
locally approximated by entire functions. This theorem cannot in
general be extended to a global result. For example, the function
$f(z,\,w)=1/z$ on the CR manifold
$\{(z,\,w)\in\mathbb{C}^2\,;\,|z|=1\}$ cannot be uniformly
approximated on compacts by entire functions.

There are, however, results indicating that global approximation
is possible in certain situations where there are no topological
obstructions. For example, Boggess and Dwilewicz
\cite{boggess-dwilewicz:hypersurface} showed that continuous CR
functions on hypersurface graphs can be approximated uniformly on
compacts by entire functions. The problem is more complicated in
higher codimension. Dwilewicz and Gauthier
\cite{dwilewicz-gauthier:approximation} have proved global
approximation results in this case. Their results require certain
convexity restrictions. Nunemacher \cite{nunemacher:approximation}
similarly proved a global approximation result, for the case of
totally real submanifolds. For rigid graphs, where the graphing
function is independent of the totally real coordinates, certain
global CR approximation results are known. (See e.g.
\cite{boggess:quadric} where approximation in $L^p$-norm is
established.) In this work, we establish global CR approximation
results on a class of graphs which contain the rigid ones, under
an additional assumption that the graphing functions satisfy a
polynomial growth condition. The class under consideration, called
the Bloom-Graham model graphs, is defined as follows.

\begin{definition}
Let $M$ be a generic CR submanifold of $\mathbb{C}^n$ of
codimension $d$, $1\leq d\leq n-1$.\\
Suppose the coordinates of $\mathbb{C}^n$ are given by
$(z,\,w)\in\mathbb{C}^d\times\mathbb{C}^{n-d}$,
\begin{alignat*}{3}
z_j&=x_j+iy_j,&\quad  j&=1,\,\ldots,\,d\\
w_j&=u_j+iv_j,&\quad  j&=1,\,\ldots,\,n-d.
\end{alignat*}
$M$ is a {\em Bloom-Graham model graph} if $M$ is given globally
as a graph of the form

\[M=\{(z,\,w)\in\mathbb{C}^n\,;\,y_j=h_j(x_1,\,\ldots,\,x_{j-1},\,w),\,j=1,\,\ldots,\,d\},\]
where
\[h=(h_1,\,\ldots,\,h_d):\mathbb{R}^d\times\mathbb{C}^{n-d}\rightarrow\mathbb{R}^d\]
is a $C^1$ map.
\end{definition}

Notice that the graphing functions are partially decoupled in
their dependence on the variables in the real directions, in the
sense that $h_j$ depends only on $x_1,\,\ldots,\,x_{j-1}$.
Additionally, these graphing functions look similar to the lower
order terms in the Bloom-Graham normal form for a CR manifold
\cite{bloom-graham:normal}. We now state our main theorem.

\begin{theorem}\label{theorem:main theorem}
Let $M$ be a Bloom-Graham model graph, as above. Assume that
\[|Dh(x,\,w)|\leq C(1+|x|^N+|w|^N),\]
for all $x\in\mathbb{R}^d$, $w\in\mathbb{C}^{n-d}$, where $D$
denotes the real derivative, and where $N$ and $C$ are uniform
constants. Let $K$ be a compact subset of $M$. Then there exists a
compact subset, $K'$, of $M$ with $K\subset K'$, such that if $f$
is a continuous CR function on a neighbourhood of $K'$, then there
is a sequence of entire functions which converge to $f$ uniformly
on $K$.
\end{theorem}

\section{Outline of Proof}
As with the proof of Baouendi and Treves' local approximation
result, our proof proceeds by using a convolution kernel which is
entire, and which is integrated along a totally real
$n$-dimensional slice of $M$.

For a point $p=(z,\,w)=(x+ih(x,\,w),\,u+iv)\in M$ define the slice

\[M_p=M_v=\{(\zeta,\,\eta)\in M\,;\,\im\eta=v\}.\]

The kernel we use is somewhat different than the Gaussian used by
Baouendi and Treves. For $\zeta\in\mathbb{C}^d$,
$\eta\in\mathbb{C}^{n-d}$ define
\[\widetilde{E}(\zeta,\,\eta)=\zeta_d^2+\sum_{j=1}^{d-1}\Lambda_j\biggl(
\zeta_j^2+\zeta_j^{P_j}\biggl)+\Gamma\biggl(\sum_{j=1}^{n-d}\eta_j^2+\eta_j^Q\biggr),\]
where the constants $\Lambda_j$, $P_j$, $\Gamma$ and $Q$ are
positive integers which will be chosen later, and which will
depend only on the given compact $K$.

For $\epsilon>0$, let
\[E_{\epsilon}(\zeta,\,\eta)=\frac{\widetilde{E}(\epsilon\zeta,\,\epsilon\eta)}{\epsilon^2}=
\zeta_d^2+\sum_{j=1}^{d-1}\Lambda_j\biggl(\zeta_j^2+\epsilon^{P_j-2}\zeta_j^{P_j}\biggr)
+\Gamma\biggl(\sum_{j=1}^{n-d}\eta_j^2+\epsilon^{Q-2}\eta_j^{Q}\biggr).\]

For $R>0$ large enough, let
$\chi_R:\mathbb{R}^d\times\mathbb{R}^{n-d}\rightarrow\mathbb{R}$
be a smooth cutoff function which is $1$ on $\{|x|+|u|\leq
R\}\supset K$ and vanishes outside of $\{|x|+|u|\leq R+1\}$.

Suppose $f$ is a continuous function defined on $M$; for\\
$p=(z,\,w)=(x+iy,\,u+iv)\in\mathbb{C}^d\times\mathbb{C}^{n-d}$
define
\begin{align*}
G_{\epsilon}(f)(z,\,w)=\frac{1}{C_1\epsilon^n}\underset{(\zeta,\,\eta)\in
M_v}{\int} &\chi_R(\re\zeta,\,\re\eta)f(\zeta,\,\eta)\\
&\exp\biggl(-E_{\epsilon}\biggl(\biggl(\frac{\zeta-z}{\epsilon}\biggr),\,\biggl(\frac{\eta-w}{\epsilon}\biggr)\biggr)
\biggr)\,d\zeta\wedge d\eta,
\end{align*}
where $C_1$ is a normalizing constant, to be chosen later.

We will show (Section \ref{section:approx}) that
$G_{\epsilon}(f)\rightarrow  f$ on $K$, as $\epsilon\rightarrow
0$. Note that the integrand defining $G_{\epsilon}(f)(z,\,w)$ is
holomorphic in $z$ and $w$. However, $G_{\epsilon}(f)(z,\,w)$ is
not necessarily holomorphic in $(z,\,w)$ since the domain of
integration, $M_v$, depends on $v=\im w$. Later (Section
\ref{section:slice}), we will show that if $f$ is CR on a suitable
subset, $K'$, of $M$ then the domain of integration $M_v$ can be
fixed, independent of $w$. The resulting sequence of functions,
denoted $F_{\epsilon}(f)(z,\,w)$, will be entire and will converge
to $f$ uniformly on $K$ as $\epsilon\rightarrow 0$.

\section{Kernel Estimates}
To prove that $G_{\epsilon}(f)\rightarrow f$, we begin by pulling
back the integral defining $G_{\epsilon}(f)$ from $M_v$ to
$\mathbb{R}^d\times\mathbb{R}^{n-d}$, via the map
$H^v:\mathbb{R}^d\times\mathbb{R}^{n-d}\rightarrow M_v$:
\begin{align*}
H^v(s,\,t)&=(\zeta^v(s,\,t),\,\eta^v(s,\,t))\\
          &=(s+ih(s,\,t+iv),\,t+iv),
\end{align*}
for $(s,\,t)\in\mathbb{R}^d\times\mathbb{R}^{n-d}$. We obtain
\begin{align*}
G_{\epsilon}(f)(z,\,w)=&
\frac{1}{C_1\epsilon^n}\underset{(s,\,t)\in\mathbb{R}^d\times\mathbb{R}^{n-d}}
{\int}\chi_R(s,\,t)f(H^v(s,\,t))\\
&\exp\biggl(-E_{\epsilon}\biggl(\biggl(\frac{\zeta^v(s,\,t)-z}{\epsilon}\biggr),\,\biggl(\frac{\eta^v(s,\,t)-w}
{\epsilon}\biggr)\biggr)\biggr)\,d s\, d t.
\end{align*}
Note that since $M$ is a Bloom-Graham graph,
$(H^{v})^*(d\zeta\wedge d\eta)=ds\wedge dt$.

We now need a key estimate on the exponent of our kernel, which we
will state in terms of $\widetilde{E}(\zeta,\,\eta)$, since
$E_{\epsilon}(\zeta/\epsilon,\,\eta/\epsilon)=\widetilde{E}(\zeta,\,\eta)/\epsilon^2$.

\begin{lemma}\label{lemma:lemma 1}
Let $M$ and $K$ be as in Theorem \ref{theorem:main theorem}. There
exist positive even integers $\Lambda_1,\,\ldots,\,\Lambda_{d-1},$
$P_1,\,\ldots,\,P_{d-1},$ $\Gamma,\,Q$, and a positive real
constant $\widetilde{C}$, all of which depend only on $K$, such
that
\begin{equation}\label{inequality:main estimate}
\begin{split}
-\re \widetilde{E}\biggl(&\zeta^{v'}(s,\,t)-z,\,
\eta^{v'}(s,\,t)-w\biggr)\\
&\leq
-\frac{(s_d-x_d)^2}{2}-\sum_{j=1}^{d-1}\biggl((s_j-x_j)^2+(s_j-x_j)^{P_j}\biggr)
\\
&\quad-\sum_{j=1}^{n-d}\biggl((t_j-u_j)^2+(t_j-u_j)^{Q}\biggr)+
\widetilde{C}\biggl(\sum_{j=1}^{n-d}\biggl((v_j'-v_j)^2+(v_j'-v_j)^Q\biggr)\biggr),\\
\end{split}
\end{equation}
for all $(z,\,w)=(x+ih(x,\,w),\,u+iv)\in K$ and all $(s,\,t+iv')$
in $\mathbb{R}^d\times\mathbb{C}^{n-d}$.
\end{lemma}

\begin{remark}
To show that $G_{\epsilon}(f)\rightarrow f$ in Section
\ref{section:approx}, we will only need the above estimate with
$v'=v$. In order to fix the domain of integration in Section
\ref{section:slice}, we will need the above estimate with $v'\neq
v$.
\end{remark}

\begin{proof}
For
$(\zeta^{v'},\,\eta^{v'})=(\zeta^{v'}(s,\,t),\,\eta^{v'}(s,\,t))$
and $(z,\,w)\in M$, we have
\begin{equation}\label{equality:coords}
\zeta^{v'}_j-z_j=(s_j-x_j)+i[h_j(s_1,\,\ldots,\,s_{j-1},\,t,\,v')-h_j(x_1,\,\ldots,\,x_{j-1},\,u,\,v)].
\end{equation}
We estimate the imaginary part by using the mean value theorem and
the assumed polynomial growth estimate on $|Dh|$ to obtain
\begin{equation}\label{inequality:mean value}
\begin{split}
|h_j(s_1,\,\ldots,\,s_{j-1},&\,t,\,v')-h_j(x_1,\,\ldots,\,x_{j-1},\,u,\,v)|\\
\leq
C\biggl(C_K&+\sum_{k=1}^{j-1}|s_k-x_k|^N+|t-u|^N+|v'-v|^N\biggr)\\
&\cdot\biggl(|t-u|+|v'-v|+\sum_{k=1}^{j-1}|s_k-x_k|\biggr),
\end{split}
\end{equation}
for all $(s,\,t+iv')\in\mathbb{R}^d\times\mathbb{C}^{n-d}$ and all
$(z,\,w)\in K$.

Note that we have used the fact that if $K\subset\mathbb{R}^d$ is
compact, then there is a constant $C_K$ such that
\[1+|s^*|^N\leq C_K(1+|s-x|^N),\]
for all $s\in\mathbb{R}^d$, $x\in K$, and $s^*\in\mathbb{R}^d$
lying between $x$ and $s$.

We will now make repeated use of a standard arithmetic inequality.
Fix any $p,\,q>1$ with $\frac{1}{p}+\frac{1}{q}=1$ and fix any
small $\delta>0$; then there is a large constant $L_{\delta}$ with
\begin{equation}\label{inequality:arithmetic}
ab\leq\delta a^p+L_{\delta}b^q,
\end{equation}
for any $a,\,b\geq 0$. A standard special case is $p=q=2$ and
$\delta=L_{\delta}=1/2$.

Using (\ref{equality:coords}) and (\ref{inequality:mean value}),
and the above arithmetic inequality, we see that if $P$ is an even
integer, then there exists a positive integer $M=M_P$, depending
only on $P$ and $N$, such that

\begin{equation}\label{inequality:zeta P}
\begin{split}
-\re\biggl(&\zeta^{v'}_j(s,\,t)-z_j\biggr)^P\leq -\frac{(s_j-x_j)^P}{2}\\
&\quad
+C_K\biggl[\sum_{k=1}^{j-1}\biggl(|s_k-x_k|^2+|s_k-x_k|^{M}\biggr)
+|t-u|^2+|t-u|^{M}\\
&\quad+|v'-v|^2+|v'-v|^{M}\biggr].
\end{split}
\end{equation}
(Inequality (\ref{inequality:arithmetic}) is used to handle cross
terms in the expansion of $(\zeta^{v'}_j-z_j)^P$.) Here and below,
$C_K$ is a constant which depends only on $K$ and may change from
line to line.

Let $j=d$ and $P=2$ in (\ref{inequality:zeta P}). The term
$-(s_d-x_d)^2/2$ on the right side of (\ref{inequality:zeta P}) is
the first term on the right side of our desired inequality
(\ref{inequality:main estimate}). Next, let $j=d-1$ with $P=2$,
and then let $j=d-1$ with $P=P_{d-1}$, an even integer which is
greater than $M=M_2$ in (\ref{inequality:zeta P}). Since the sum
on the right hand side of (\ref{inequality:zeta P}) does not
involve $(s_l-x_l)$ for $l\geq j$, we can choose a constant
$\Lambda_{d-1}>2(C_K+1)$ and then combine these three inequalities
(i.e. (\ref{inequality:zeta P}) with $j=d$, $P=2$; $j=d-1$, $P=2$;
and $j=d-1$ and $P=P_{d-1}>M_2$) to obtain
\begin{equation}\label{inequality:lemma estimate}
\begin{split}
-\re&\biggl[\biggl(\zeta_d^{v'}(s,\,t)-z_d\biggr)^2+\Lambda_{d-1}
\biggl[\biggl(\zeta_{d-1}^{v'}(s,\,t)-z_{d-1}\biggr)^2\\
&\quad\quad\quad+\biggl(\zeta_{d-1}^{v'}(s,\,t)-z_{d-1}\biggr)^{P_{d-1}}
\biggr]\biggr]\\
&\leq -\frac{(s_d-x_d)^2}{2}-(s_{d-1}-x_{d-1})^2-(s_{d-1}-x_{d-1})^{P_{d-1}}\\
&\quad\quad\quad+C_K\biggl[\sum_{k=1}^{d-2}\biggl(|s_k-x_k|^2+|s_k-x_k|^{\widehat{M}}\biggr)
+|t-u|^2+|t-u|^{\widehat{M}}+|v'-v|^2\\
&\quad\quad\quad+|v'-v|^{\widehat{M}}\biggr],
\end{split}
\end{equation}
where $\widehat{M}$ depends only on $P_{d-1}$ and $M=M_2$, which
in turn depend only on $N$. Note that the first three terms on the
right of (\ref{inequality:lemma estimate}) agree with the first
three terms on the right of (\ref{inequality:main estimate}).
Continuing in this manner, we can inductively choose constants
$\Lambda_{d-1},\,P_{d-1},\,\ldots,\,\Lambda_1$, $P_{1}$, and $M'$
an even integer, in that order and depending only on $K$, so that

\begin{equation}\label{inequality:zeta terms}
\begin{split}
-\re&\biggl[\biggl(\zeta_d^{v'}(s,\,t)-z_d\biggr)^2
+\sum_{j=1}^{d-1}\Lambda_{j}
\biggl[\biggl(\zeta_{j}^{v'}(s,\,t)-z_{j}\biggr)^2+\biggl(\zeta_{j}^{v'}(s,\,t)-z_{j}\biggr)^{P_{j}}
\biggr]\biggr]\\
&\leq -(s_d-x_d)^2-\biggl(\sum_{j=1}^{d-1}(s_{j}-x_{j})^2+(s_{j}-x_{j})^{P_{j}}\biggr)\\
&+C_K\biggl[|t-u|^2+|t-u|^{M'}+|v'-v|^2+|v'-v|^{M'}\biggr].
\end{split}
\end{equation}

Now write $w=u+iv$, $\eta^{v'}=t+iv'$. Using the arithmetic
inequality (\ref{inequality:arithmetic}), and with $Q$ an even
integer, we obtain
\begin{equation*}
-\re(\eta^{v'}(s,\,t)-w)^Q\leq \frac{-(t-u)^Q}{2}+L(v'-v)^Q,
\end{equation*}
where $L$ is a fixed constant. Using (\ref{inequality:zeta
terms}), and choosing $\Gamma>2(C_K+1)$ and $Q\geq M'$, we obtain
Lemma \ref{lemma:lemma 1}.
\end{proof}

\section{Approximation to the Identity}\label{section:approx}
\begin{lemma}\label{lemma:approximation} If $f$ is a continuous function on $M$, then
$G_{\epsilon}(f)\rightarrow f$ on $K$ as $\epsilon\rightarrow 0$.
\end{lemma}

\begin{proof}
Recall that the domain of integration of $G_{\epsilon}(f)(z,\,w)$
is $M_v=M\cap\{v'=\im\eta=v\}$ where $v=\im w$.

After making the change of variables
\begin{equation*}
s=x-\epsilon\hat{s},\quad t=u-\epsilon\hat{t},
\end{equation*}
we obtain
\begin{align*}
G_{\epsilon}&(f)(z,\,w)=\frac{1}{C_1}\underset{(\hat{s},\,\hat{t})\in
\mathbb{R}^d\times\mathbb{R}^{n-d}}{\int}
\chi_R(x-\epsilon\hat{s},\,u-\epsilon\hat{t})f(H^v(x-\epsilon\hat{s},\,u-\epsilon\hat{t}))\\
&\exp\biggl(-E_{\epsilon}\biggl(\biggl(\frac{\zeta^v(x-\epsilon\hat{s},\,u-\epsilon\hat{t})-z}{\epsilon}\biggr),
\,\biggl(\frac{\eta^v(x-\epsilon\hat{s},\,u-\epsilon\hat{t})-w}
{\epsilon}\biggr)\biggr)\biggr)\,d \hat{s}\,d \hat{t}.
\end{align*}

Using Lemma \ref{lemma:lemma 1} with $v'=v$, and the fact that
$E_{\epsilon}(\zeta/\epsilon,\,\eta/\epsilon)=\widetilde{E}(\zeta,\,\eta)/\epsilon^2$,
we see that the real part of the exponent of our kernel is less
than or equal to
\begin{equation*}
-\frac{\hat{s}_d^2}{2}-\sum_{j=1}^{d-1}\biggl(\hat{s}_j^2+\epsilon^{P_j-2}\hat{s}_j^{P_j}\biggr)
-\sum_{j=1}^{n-d}\biggl(\hat{t}_j^2+\epsilon^{Q-2}\hat{t}_j^{Q}\biggr)
\leq -\frac{\hat{s}_d^2}{2}-\sum_{j=1}^{d-1}\hat{s}_j^2
-\sum_{j=1}^{n-d}\hat{t}_j^2
\end{equation*}

Since $\chi_R\cdot f$ is bounded, the Dominated Convergence
Theorem allows us to let $\epsilon\rightarrow 0$ in the integrand
of $G_{\epsilon}$. The resulting integral is
\begin{equation}\label{equation:constant}
\begin{split}
\frac{1}{C_1}&\underset{(\hat{s},\,\hat{t})\in
\mathbb{R}^d\times\mathbb{R}^{n-d}}{\int}
\chi_R(x,\,u)f(H^v(x,\,u)) \exp [-E_0(D
H^v(x,\,u)\cdot(\hat{s},\,\hat{t})^T)] \,d \hat{s}\,d \hat{t},
\end{split}
\end{equation}
where
\[E_0=\zeta_d^2+\sum_{j=1}^{d-1}\Lambda_j\zeta_j^2
+\sum_{j=1}^{n-d}\Gamma_j\eta_j^2,\] and where $D$ is the usual
first order derivative with respect to $x$ and $u$.

We must show that
\begin{equation}\label{equation:complex integral}
\begin{split}
&\underset{(\hat{s},\,\hat{t})\in
\mathbb{R}^d\times\mathbb{R}^{n-d}}{\int} \exp[-E_0(D
H^v(x,\,u)\cdot(\hat{s},\,\hat{t})^T)] \,d \hat{s}\,d \hat{t}
\end{split}
\end{equation}
is a constant, independent of $x$, $u$, and $v$. For then by
letting $C_1$ be this constant, the expression in
(\ref{equation:constant}) becomes $f(H^v(x,\,u))=\chi_R
f(z,\,w)=f(z,\,w)$ for $(z,\,w)\in K$, and the proof of Lemma
\ref{lemma:approximation} will be complete.

Note that since $M$ is a Bloom-Graham model graph, $DH^v$ is of
the form \begin{equation*} \mathcal{M}(A,\,B)=
\begin{pmatrix}
I_{d\times d} &0\\
0 & I_{(n-d)\times (n-d)}\\
\end{pmatrix}
+
\begin{pmatrix}
A& B \\
0 & 0\\
\end{pmatrix},
\end{equation*}
where $I_j$ is the $j\times j$ identity matrix, $A$ is a $d\times
d$ lower triangular matrix with zeros along the diagonal, and $B$
is a $d\times (n-d)$ matrix. In particular,
$\det\mathcal{M}(A,\,B)=1$. By increasing the values of
$\Lambda_j$ and $\Gamma_j$ if necessary (as in the proof of Lemma
\ref{lemma:lemma 1}, with only quadratic terms in the exponent),
we see that
\begin{equation}\label{equation:integral}
\begin{split}
&\underset{(\hat{s},\,\hat{t})\in
\mathbb{R}^d\times\mathbb{R}^{n-d}}{\int}
\exp[-E_0(\mathcal{M}(A,\,B)\cdot(\hat{s},\,\hat{t})^T)] \,d
\hat{s}\,d \hat{t}
\end{split}
\end{equation}
is complex analytic in the entries of $A$ and $B$ on a complex
neighbourhood of the region given by
\[|A_{j,\,k}|\leq 2\max_{
K}\biggl(\biggl|\frac{\partial h_j(x,\,u,\,v)}{\partial
x_k}\biggr|\biggr), \quad 1\leq j<k\leq d\] and
\[|B_{j,\,k}|\leq 2\max_{
K}\biggl(\biggl|\frac{\partial h_j(x,\,u,\,v)}{\partial
u_k}\biggr|\biggr),\quad 1\leq j\leq d,\,1\leq k\leq n-d.\]

If $A$ and $B$ are real-valued, a change of variables of the form
$(s,\,t)^T=\mathcal{M}(A,\,B)\cdot(\hat{s},\,\hat{t})^T$ shows
that the integral in (\ref{equation:integral}) is independent of
$A$ and $B$ (again, recall that $\det\mathcal{M}(A,\,B)=1$). By
the identity theorem for holomorphic functions, the same is true
when $A$ and $B$ are complex valued, belonging to the above
neighbourhood. Thus (\ref{equation:complex integral}) is
independent of $x$, $u$, and $v$, as desired.
\end{proof}

\section{Fixing the Slice}\label{section:slice} We have established that
$G_{\epsilon}(f)(z,\,w)\rightarrow f(z,\,w)$ as
$\epsilon\rightarrow 0$ for $(z,\,w)\in K$. While the integrand in
$G_{\epsilon}$ is holomorphic in $z$ and $w$, the domain of
integration depends on $v=\im w$. Thus $G_{\epsilon}$ is not
necessarily a holomorphic function. To remedy this defect, we fix
the domain of integration at $M_{v_0}$, independent of $v=\im w$,
and define
\begin{align*}
F_{\epsilon}(f)(z,\,w)=&\frac{1}{C_1\epsilon^n}\underset{(\zeta,\,\eta)\in
M_{v_0}}{\int} \chi_R(\zeta,\,\eta)f(\zeta,\,\eta)\\
&\quad\quad\exp\biggl(-E_{\epsilon}\biggl(\biggl(\frac{\zeta-z}{\epsilon}\biggr),
\,\biggl(\frac{\eta-w}{\epsilon}\biggr)\biggr)
\biggr)\,d\zeta\wedge d\eta,
\end{align*}
with $v_0$ a fixed point in the projection of the compact $K$ onto
the $v$-axis. As $\chi_R\cdot f$ is compactly supported, this
integral is well defined, and thus holomorphic, for all
$(z,\,w)\in\mathbb{C}^n$, $\epsilon>0$.

\begin{lemma} Let $K$ be a compact
subset of $M$. Let $R'$ be any number larger than
$\max_{(z,\,w)\in K}|v|$. Then there exists $R>0$ and $C>0$ such
that if $f$ is a continuous CR function on a neighbourhood of
$K'=M\cap[\{|x|+|u|\leq R+1\}\times\{|v|<R'\}],$ then
\[|F_{\epsilon}(f)(z,\,w)-G_{\epsilon}(f)(z,\,w)|\leq C\epsilon,\]
for all $(z,\,w)\in K$, and all $\epsilon>0$.
\end{lemma}

Since $F_{\epsilon}(f)$ is entire and $G_{\epsilon}(f)\rightarrow
f$ on $K$, clearly this lemma completes the proof of Theorem
\ref{theorem:main theorem}.

\begin{proof}
Consider the manifold $\widetilde{M}_v$ defined as
\[\widetilde{M}_v=\{(\zeta,\,\eta)\in M\,;\,\im\eta=v'=rv+(1-r)v_0,\,0\leq r \leq1\}.\]
$\widetilde{M}_v$ is an $(n+1)$ real dimensional submanifold of
$M$ with boundary components $M_v$ and $M_{v_0}$. By Stokes'
theorem
\begin{equation*}
\begin{split}
G_{\epsilon}(f)(z,\,w) &=F_{\epsilon}(f)(z,\,w)\\
&+\frac{1}{C_1\epsilon^n}\underset{(\zeta,\,\eta)\in
\widetilde{M}_v}{\int}d_{(\zeta,\,\eta)}\biggl[
\chi_R(\zeta,\,\eta)f(\zeta,\,\eta) \\
&\quad\quad\quad\quad
\exp\biggl(-E_{\epsilon}\biggl(\biggl(\frac{\zeta-z}{\epsilon}\biggr),\,\biggl(\frac{\eta-w}{\epsilon}\biggr)\biggr)
\biggr)\,d\zeta\wedge d\eta\biggr].\\
\end{split}
\end{equation*}
The presence of $d\zeta\wedge d\eta$ implies that only
$\overline{\partial}_{(\zeta,\,\eta)}$ terms appear. If $f$ is CR
on a neighbourhood of $K'$, then the support of the integrand is
contained in $\{(\zeta,\,\eta)=(s+ih(s,\,t,\,v),\,t+iv)\in
M\,;\,R\leq |s|+|t|\leq R+1,\,|v|<R'\}$.

Since $|x|$, $|u|$ and $|v'-v|$ (with $v'=rv+(1-r)v_0$) are
bounded in terms of the diameter of $K$, inequality
(\ref{inequality:main estimate}) shows that choosing $R$ suitably
large relative to the diameter of $K$ ensures that
\[\exp\biggl(-E_{\epsilon}\biggl(\biggl(\frac{\zeta-z}{\epsilon}\biggr),\,\biggl(\frac{\eta-w}{\epsilon}\biggr)\biggr)
\biggr)\leq e^{\frac{-R^2}{4\epsilon^2}},\] for $(\zeta,\,\eta)$
in the support of $\overline{\partial}\chi_R$ and $(z,\,w)\in K$.
\end{proof}

\begin{remark}
Though the assumption of $M$ being a Bloom-Graham model graph is
referenced throughout this work, this assumption is only
critically used in the proof of Lemma \ref{lemma:lemma 1},
specifically in (\ref{inequality:zeta P}) and
(\ref{inequality:lemma estimate}).
\end{remark}
\bibliographystyle{amsplain}
\bibliography{bibliography}

\end{document}